\documentclass[11pt]{article}
\usepackage{amsmath, amsthm, amssymb}
\usepackage{enumerate}
\usepackage[top=30truemm, bottom=30truemm, left=25truemm, right=25truemm]{geometry}
\usepackage[dvipdfm]{graphicx}
\usepackage{accents}
\usepackage{bm}
\usepackage{caption}
\usepackage{cases}
\usepackage{tikz-cd}
\usepackage{titlesec}
\usepackage[all]{xy}

\makeatletter
\def\sbar{\accentset{{\cc@style\underline{\mskip9mu}}}}
\def\mbar{\accentset{{\cc@style\underline{\mskip12mu}}}}
\def\lbar{\accentset{{\cc@style\underline{\mskip18mu}}}}

\makeatother
\newcommand{\keywords}[1]{\textbf{{Keywords: }} #1}
\newcommand{\MSC}[1]{\textbf{{2010 Mathematical Subject Classification: }} #1}
\titleformat*{\section}{\normalfont\large\bfseries}

\theoremstyle{definition}
\newtheorem{Def}{Definition}[section]
\newtheorem{Thm}[Def]{Theorem}
\newtheorem{Prop}[Def]{Proposition}
\newtheorem{Lem}[Def]{Lemma}
\newtheorem{Cor}[Def]{Corollary}
\newtheorem{Rmk}[Def]{Remark}
\newtheorem{Exp}[Def]{Example}
\title{\vspace{-5mm}On the maximality of genus-3 nonhyperelliptic curves of Ciani type}
\author{Ryo Ohashi}

\begin{document}
\maketitle\vspace{-6mm}
\begin{abstract}
In this paper, we study a Ciani curve $C: x^4 + y^4 + z^4 + rx^2y^2 + sy^2z^2 + tz^2x^2 = 0$ in positive characteristic $p \geq 3$. We will show that if $C$ is superspecial, then its standard form is maximal or minimal over $\mathbb{F}_{p^2}\hspace{-0.3mm}$ without taking its $\mathbb{F}_{p^2}$-form.
\end{abstract}

\keywords{Algebraic curve, Superspecial curve, Maximal curve, Positive characteristic}\par
\MSC{14G05, 14G17, 14H45, 14H50}

\section{Introductrion}
Throughout this paper, a curve always means a projective variety in positive characteristic $p \geq 3$ of dimension one. It is well-known that all nonsingular genus-$g$ curves $C$ defined over $\mathbb{F}_q$ with $q = p^n$ satisfy the Hasse-Weil inequality\vspace{-1mm}
\[
	1+ q - 2g\sqrt{q} \leq \#C(\mathbb{F}_q) \leq 1 + q  + 2g\sqrt{q},\vspace{-1mm}
\]
where $C(\mathbb{F}_q)$ denotes the set of $\mathbb{F}_q$-rational points of $C$. Now, a nonsingular curve $C$ is called \textit{maximal} (resp.\ \textit{minimal}) when the number of $\mathbb{F}_q$-rational points of $C$ attains the upper (resp.\ lower) bound. Maximal curves have been investigated for their applications to coding theory. On the other hand, we call a curve $C$ \textit{superspecial} if ${\rm Jac}(C)$ is isomorphic to the product of supersingular elliptic curves. This is equivalent to saying that the $a$-number of $C$ is equal to $g$. The $a$-number of $C$ is defined to be the dimension of ${\rm Hom}(\alpha_p,{\rm Jac}(C)[p])$ where $\alpha_p$ is the kernel of the Frobenius map on the additive group $\mathbb{G}_a$. And then, it is also known that any maximal or minimal curve over $\mathbb{F}_{p^2}$ is superspecial, though a superspecial curve over $\mathbb{F}_{p^2}$ is not necessarily maximal nor minimal.\par

In this paper, we study the maximality of genus-3 nonhyperelliptic curves of Ciani type, and we will call them Ciani curves briefly. A Ciani curve $C$ is a plane quartic defined by the equation
\[
	C: x^4 + y^4 + z^4 + rx^2y^2 + sy^2z^2 + tz^2x^2 = 0,
\]
which was studied by Ciani \cite{Ciani}. Brock \cite[Theorem 3.15]{Brock} studied the superspeciality of Ciani curves and their enumerations, using the result by Hashimoto \cite{Hashimoto} on the computations of the class numbers of quaternion unitary groups. Our main result is as below:\vspace{-1mm}
\begin{Thm}
Assume that a nonsingular curve 
\[
	C: x^4 + y^4 + z^4 + rx^2y^2 + sy^2z^2 + tz^2x^2 = 0
\]
is superspecial, then $r,s$ and $t$ belong to $\mathbb{F}_{p^2}$. Moreover $C$ is maximal or minimal over $\mathbb{F}_{p^2}$.
\end{Thm}\vspace{-1mm}
\noindent See Corollary 4.6 for a condition determining whether $C$ is maximal or minimal.

The remainder of this paper is structured as follows. In Section 2, we review several properties of Ciani curves. Proposition 2.5 gives us the classification of automorphism groups of a Ciani curve. In Section 3, we look into the structure of a Ciani curve $C$. In particular, we describe explicitly the elliptic curves appearing as quotients of $C$ by involutions. In Section 4, we prove Theorem 1.1.\vspace{-1mm}

\subsection*{Acknowledgments}\vspace{-1mm}
This paper was written while the author is a Ph.D.\,student at Yokohama National University, and I would like to express my special thanks to my supervisor Prof.\ Shushi Harashita for his guidance.

\section{Ciani curve}
\setcounter{equation}{0}
Let $K$ be a perfect field of characteristic $p \geq 3$. In this section, we consider a nonhyperelliptic curve of genus 3 defined over $K$\vspace{-1mm}
\begin{equation}
	C: x^4 + y^4 + z^4 + rx^2y^2 + sy^2z^2 + tz^2x^2 = 0,
\end{equation}
which is called a \textit{Ciani curve}. First of all, we discuss the singularity of a Ciani curve.
\begin{Lem}
The curve $C$ is nonsingular if and only if $r,s,t \neq \pm 2$ and $r^2+s^2+t^2-rst-4 \neq 0$.
\begin{proof}
Put $F := x^4 + y^4 + z^4 + rx^2y^2 + sy^2z^2 + tz^2x^2$. Then, we have
\[
	\frac{\partial{F}}{\partial{x}} = 2x(2x^2 + ry^2 + tz^2), \quad \frac{\partial{F}}{\partial{y}} = 2y(2y^2 + rx^2 + sz^2), \quad \frac{\partial{F}}{\partial{z}} = 2z(2z^2 + sy^2 + tx^2).\vspace{-1mm}
\]
To show the ``if''-part, assume that $(x:y:z) \in \mathbb{P}^2$ is a singular point of $C$.\par
Firstly, we note that there is no point $(x:y:z)$ on $C$ such that two among $x,y$ and $z$ are zero; e.g., if $x = y = 0$, it follows from $F=0$ that $z=0$, which is a contradiction.\par
Secondly, consider a singular point $(x:y:z)$ on $C$ such that only one among $x,y$ and $z$ is zero, e.g.; the case $x = 0,\,y \neq 0$ and $z \neq 0$. It follows from $2y^2 + sz^2 = 2z^2 + sy^2 = 0$ that $(4-s^2)z^2 = 0$, whence $s = \pm 2$. One can check the case of $y=0$ or $z=0$ in the same way.\par
Lastly, consider a singular point $(x:y:z)$ on $C$ such that all of $x,y$ and $z$ are non-zero. Using the Jacobian criterion, we obtain\vspace{-1mm}
\[
\left( \begin{array}{ccc}
	2 & r & t\\
	r & 2 & s \\
	t & s & 2 \\ 
	\end{array} \right)
	\left( \begin{array}{c}
	x^2 \\
	y^2 \\
	z^2 \\ 
	\end{array} \right)
    = \left( \begin{array}{c}
	0 \\
	0 \\
	0 \\ 
	\end{array} \right).\vspace{-1mm}
\]
Here, we can calculate\vspace{-2mm}
\[
	\det\left( \begin{array}{ccc}
	2 & r & t\\
	r & 2 & s \\
	t & s & 2 \\ 
	\end{array} \right) = -2(r^2+s^2+t^2 - rst - 4),
\]
and thus $r^2+s^2+t^2 - rst - 4 = 0$ since $(x^2,y^2,z^2) \neq (0,0,0)$.\par
Conversely, we show the ``only if''-part. If $r = \pm2$, we can find singular points $(x:y:0)$ on $C$ such that $x^2 \pm y^2 = 0$. One can check the case of $s = \pm 2$ or $t = \pm 2$ in the same way. Moreover, by a tedious calculation, we can find singular points $(x:y:z)$ on $C$ such that
\[
	(4-r^2)(rt-2s)x^2 = (4-r^2)(rs-2t)y^2 = (rt-2s)(rs-2t)z^2
\]
if $r^2+s^2+t^2-rst-4 = 0$. The proof is completed.
\end{proof}
\end{Lem}
\begin{Def}
We say that a nonsingular Ciani curve of the form (2.1) is of \textit{$[a,b,c]$-type} when the triple $(r,s,t)$ is a permutation of
\[
	(a,b,c),(a,-b,-c),(-a,b,-c),(-a,-b,c).
\]
If both $C$ and $C'$ are two Ciani curves of $[a,b,c]$-type, it is obvious that $C$ is isomorphic to $C'$.
\end{Def}
\newpage
\begin{Rmk}
Suppose that $C$ and $C'$ are two Ciani curves, then the type of $C$ is not necessarily same as that of $C'$ even if $C$ and $C'$ are isomorphic. For example, the Fermat curve $x^4+y^4+z^4 = 0$ is isomorphic to the curve $x^4 + y^4 + 6y^2z^2 + z^4 = 0$. We find a necessary and sufficient condition for two Ciani curves $C$ and $C'$ to be isomorphic in \cite[Proposition 2.2]{Howe}, but do not use it in this paper.
\end{Rmk}
\begin{Exp}
Consider a nonsingular curve
\begin{equation}
	C: x^4 + y^4 + z^4 + rx^2yz +sy^2z^2 = 0.
\end{equation}
Note that the automorphism group of $C$ contains ${\rm D}_8$ as shown in \cite[Theorem 3.1]{Ritzenthaler}. Let us confirm that this curve $C$ is of $[a,b,a]$-type with $a=\frac{r}{\sqrt{s+2}}$ and $b = 2-\frac{16}{s+2}$.\par
By replacing $u = (y+z)/2$ and $z = (y-z)/(2\sqrt{-1})$, we obtain $y = u+\sqrt{-1}v$ and $z = u+\sqrt{-1}v$ thus\vspace{-2mm}
\begin{align*}
	yz &= (u+\sqrt{-1}v)(u+\sqrt{-1}v) = u^2 + v^2,\\
	y^2 + z^2 &= (u+\sqrt{-1}v)^2 + (u+\sqrt{-1}v)^2 = 2(u^2-v^2).
\end{align*}
The left side of (3.2) can be transformed into
\begin{align*}
	x^4 + y^4 + z^4 + rx^2yz + sy^2z^2 &= x^4 + (y^2 + z^2)^2 + rx^2yz + (s-2)y^2z^2\\
	&= x^4 + (s+2)u^4 + (s+2)v^4 + rx^2u^2 + 2(s-6)u^2v^2 + rv^2x^2,
\end{align*}
and the map
\[
	u \mapsto \frac{u}{\sqrt[4]{s+2}} ,\quad v \mapsto \frac{v}{\sqrt[4]{s+2}}
\]
transforms the curve $C$ into
\[
	x^4 + u^4 + v^4 + \frac{r}{\sqrt{s+2}}x^2u^2 + \Bigl(2-\frac{16}{s+2}\Bigr)u^2v^2 + \frac{r}{\sqrt{s+2}}v^2x^2.
\]
Hence, the curve $C$ is a Ciani curve of $[a,b,a]$-type.
\end{Exp}
Next result tells us the classification of automorphism groups of Ciani curves. Let ${\rm C}_n$\hspace{1mm}(resp.\ ${\rm D}_n$) be the cyclic (resp.\ dihedral) group of order $n$ and ${\rm S}_n$ the symmetric group of degree $n$.
\begin{Prop}
The automorphism group of a Ciani curve $C$ is either of the following 7 groups:
\begin{enumerate}
\item ${\rm D}_4$. In case $C$ can be written as $[a,b,c]$-type for some $a,b,c \in \mbar{K}$.
\item ${\rm D}_8$. In case $C$ can be written as $[a,b,a]$-type for some $a,b \in \mbar{K}$.
\item ${\rm G}_{16}$. In case $C$ can be written as $[0,b,0]$-type for some $b \in \mbar{K}$.
\item ${\rm S}_4$. In case $C$ can be written as $[a,a,a]$-type for some $a \in \mbar{K}$.
\item ${\rm G}_{48}$. In case $C$ can be written as $[0,2\sqrt{-3},0]$-type.
\item ${\rm G}_{96}$. In case $C$ can be written as $[0,0,0]$-type.\vspace{-1mm}
\item ${\rm G}_{168}$. In case $C$ can be written as $\bigl[\frac{-3+\sqrt{-63}}{2},\frac{-3+\sqrt{-63}}{2},\frac{-3+\sqrt{-63}}{2}\bigr]$-type.
\end{enumerate}
\begin{proof}
It follows from \cite[Theorem 3.1]{Ritzenthaler} about (1),\,(3),\,(4) and (6). It is well-known \cite[Section 6-8]{Meagher} about (5) and (7). It is clear about (2) by Example 2.4.
\end{proof}
\end{Prop}

\section{Elliptic curves associated to a Ciani curve}
\setcounter{equation}{0}
Let $\sigma_1,\sigma_2,\sigma_3$ be automorphisms on a nonsingular Ciani curve $C$ as below:
\begin{align*}
	\sigma_1&: (x:y:z) \mapsto (-x:y:z),\\
	\sigma_2&: (x:y:z) \mapsto (x:-y:z),\\
	\sigma_3&: (x:y:z) \mapsto (x:y:-z).
\end{align*}
Put $E_i := C/\langle\sigma_i\rangle$ for $i \in \{1,2,3\}$, then we obtain the equations
\begin{align*}
	E_1 &: X^2 + y^4 + z^4 + rXy^2 + sy^2z^2 + tz^2X = 0,\\
	E_2 &: x^4 + Y^2 + z^4 + rx^2Y + sYz^2 + tz^2x^2 = 0,\\
	E_3 &: x^4 + y^4 + Z^2 + rx^2y^2 + sy^2Z + tZx^2 = 0
\end{align*}
with $X = x^2,\,Y = y^2$ and $Z = z^2$. One can easily check each $E_i$ is a genus-1 curve, and the next lemma gives the Legendre forms of $E_i$. We choose $\alpha,\beta$ and $\gamma$ such that 
\begin{equation}
	\alpha^2 = r^2-4, \quad \beta^2 = s^2-4, \quad \gamma^2 = t^2 - 4.
\end{equation}
once and fix them throughout this paper.
\begin{Lem}
We can transform $E_i$ into the Legendre forms $y^2 = x(x-1)(x-\lambda_i)$ with
\begin{align}
	\lambda_1 &= \frac{(rt-2s)-\gamma\alpha}{(rt-2s)+\gamma\alpha},\nonumber\\
	\lambda_2 &= \frac{(sr-2t)-\alpha\beta}{(sr-2t)+\alpha\beta},\\
	\lambda_3 &= \frac{(ts-2r)-\beta\gamma}{(ts-2r)+\beta\gamma}.\nonumber
\end{align}
\begin{proof}
In this proof, we will show about $E_1$. Firstly, we get the equation
\begin{align*}
	\Bigl(X + \frac{r}{2}y^2 + \frac{t}{2}z^2\Bigr)^{\!2} &= \Bigl(\frac{r}{2}y^2 + \frac{t}{2}z^2\Bigr)^{\!2} - y^4 - sy^2z^2 - z^4\\
	&= \Bigl(\frac{r^2}{4}-1\Bigr)y^4 + \Bigl(\frac{rt}{2} - s\Bigr)y^2z^2 + \Bigl(\frac{t^2}{4}-1\Bigr)z^4.
\end{align*}
By replacing $u = y$ and $v = X + \frac{r}{2}y^2 + \frac{t}{2}z^2$ and multiplying $(r^2-4)/4$ of both sides, we obtain
\[
	v^2 = u^4 + \frac{2rt-4s}{r^2-4}u^2z^2 +\frac{t^2-4}{r^2-4}z^4.
\]
Here, the right side is transformed into
\begin{align*}
	u^4+ \frac{2rt-4s}{r^2-4}u^2z^2 + \frac{t^2-4}{r^2-4}z^4 &= \Bigl(u^2+\frac{\gamma}{\alpha}z^2\Bigr)^{\!2} \!- 2\frac{\gamma\alpha - (rt-2s)}{r^2-4}u^2z^2\\
	&=  \Bigl(u^2+ buz + \frac{\gamma}{\alpha}z^2\Bigr)\Bigl(u^2- buz + \frac{\gamma}{\alpha}z^2\Bigr)
\end{align*}
with $b^2 = 2 \cdot \frac{\gamma\alpha - (rt-2s)}{r^2-4}$. Put\vspace{-2mm}
\[
	d_1 = \frac{b + \sqrt{b^2 - 4 \cdot \frac{\gamma}{\alpha}}}{2}, \quad d_2 = \frac{b - \sqrt{b^2 - 4 \cdot \frac{\gamma}{\alpha}}}{2},
\]
then we have the factorization
\[
	\Bigl(u^2+ buz + \frac{\gamma}{\alpha}z^2\Bigr)\Bigl(u^2- buz + \frac{\gamma}{\alpha}z^2\Bigr) = (u+ d_1z)(u+d_2z)(u-d_1z)(u-d_2z).
\]
The map
\[
	u \mapsto \frac{u-d_1z}{u+d_1z} \cdot \frac{d_2+d_1}{d_2-d_1}
\]
transforms this elliptic curve into the Legendre form:
\[
	v^2 = u(u-1)\biggl(u-\frac{(d_2+d_1)^2}{(d_2-d_1)^2}\biggr) = u(u-1)\biggl(u-\frac{b^2}{b^2 - 4 \cdot \frac{\gamma}{\alpha}}\biggr) = u(u-1)(u-\lambda_1).
\]
This is the desired conclusion.
\end{proof}
\end{Lem}

Here, we can regard the curve
\[
	P: X^2 + Y^2 + Z^2 + rXY + sYZ + tZX = 0
\]
as the quotient $C/\langle\sigma_1,\sigma_2\rangle$. One can confirm that $P$ is nonsingular since $C$ is nonsingular and the genus of $P$ is $0$. Hence, the curve $P$ is isomorphic to the projective line $\mathbb{P}^1$. For the above discussion, we obtain the following diagram:\vspace{-3mm}
\begin{equation}
	\vcenter{
	\xymatrix{
	& C \ar[ld] \ar[d] \ar[rd] &\\
	E_1 \ar[rd] & E_3 \ar[d] & E_2 \ar[ld]\\
	& P &
	}}\vspace{-1mm}
\end{equation}
The diagram (3.3) induces an isogeny ${\rm Jac}(C) \rightarrow E_1 \hspace{-0.3mm}\times\hspace{-0.3mm} E_2 \hspace{-0.3mm}\times\hspace{-0.3mm} E_3$ by \cite[Section 3]{Kani} of degree $2^3$. Since the degree of the isogeny is not divided by $p$, then we obtain ${\rm Jac}(C)[p] \cong (E_1 \hspace{-0.3mm}\times\hspace{-0.3mm} E_2 \hspace{-0.3mm}\times\hspace{-0.3mm} E_3)[p]$. Hence, a Ciani curve $C$ is superspecial if and only if $E_1,E_2$ and $E_3$ are supersingular.

\begin{Prop}
The reverse transformation of that in Lemma 3.1 is given by
\begin{align*}
	r &= \frac{\lambda_1\lambda_2 - \lambda_2\lambda_3 - \lambda_3\lambda_1 + 1}{\sqrt{\lambda_1\lambda_2}(1-\lambda_3)},\\
	s &=\frac{\lambda_2\lambda_3 - \lambda_3\lambda_1 - \lambda_1\lambda_2 + 1}{\sqrt{\lambda_2\lambda_3}(1-\lambda_1)},\\
	t &= \frac{\lambda_3\lambda_1 - \lambda_1\lambda_2 - \lambda_2\lambda_3 + 1}{\sqrt{\lambda_3\lambda_1}(1-\lambda_2)}.
\end{align*}
\begin{proof}
In this proof, we will show the case of $\lambda_1,\lambda_2,\lambda_3 \neq -1$. By linear fractional transformations, one can check that\vspace{-1mm}
\begin{equation}
	\frac{1+\lambda_1}{1-\lambda_1} = \frac{rt-2s}{\gamma\alpha},\quad \frac{1+\lambda_2}{1-\lambda_2} = \frac{sr-2t}{\alpha\beta}, \quad \frac{1+\lambda_3}{1-\lambda_3} = \frac{ts-2r}{\beta\gamma}.\nonumber
\end{equation}
Hence we obtain the equation
\[
	\frac{(1-\lambda_1)(1-\lambda_2)(1+\lambda_3)}{(1+\lambda_1)(1+\lambda_2)(1-\lambda_3)} = \frac{(r^2-4)(ts-2r)}{(rt-2s)(sr-2t)}.
\]
Since $(r^2-4)(ts-2r) + (rt-2s)(sr-2t) = -2r(r^2+s^2+t^2-rst-4)$, then we have
\begin{equation}
	1 + \frac{(1-\lambda_1)(1-\lambda_2)(1+\lambda_3)}{(1+\lambda_1)(1+\lambda_2)(1-\lambda_3)} =-2r \cdot \frac{r^2+s^2+t^2-rst-4}{(rt-2s)(sr-2t)}.
\end{equation}

On the other hand, note that
\begin{align}
	(rt-2s)^2 -4(r^2+s^2+t^2-rst-4) &= (r^2-4)(t^2-4),\nonumber\\
	(sr-2t)^2 -4(r^2+s^2+t^2-rst-4) &= (s^2-4)(r^2-4),\\
	(ts-2r)^2 -4(r^2+s^2+t^2-rst-4) &= (t^2-4)(s^2-4),\nonumber
\end{align}
thus one can check that\vspace{-3mm}
\begin{align*}
	\frac{\lambda_1}{(1+\lambda_1)^2} &= \frac{r^2+s^2+t^2-rst-4}{(rt-2s)^2},\\
	\frac{\lambda_2}{(1+\lambda_2)^2} &= \frac{r^2+s^2+t^2-rst-4}{(sr-2t)^2},\\
	\frac{\lambda_3}{(1+\lambda_3)^2} &= \frac{r^2+s^2+t^2-rst-4}{(ts-2r)^2}.
\end{align*}
Therefore, we have\vspace{-1mm}
\begin{equation}
	\frac{\sqrt{\lambda_1\lambda_2}}{(1+\lambda_1)(1+\lambda_2)} = \frac{r^2+s^2+t^2-rst-4}{(rt-2s)(sr-2t)}.
\end{equation}
Using the equations (3.4) and (3.6), then\vspace{-1mm}
\[
	\frac{(1+\lambda_1)(1+\lambda_2)(1-\lambda_3) + (1-\lambda_1)(1-\lambda_2)(1+\lambda_3)}{(1+\lambda_1)(1+\lambda_2)(1-\lambda_3)} = -2r \cdot \frac{\sqrt{\lambda_1\lambda_2}}{(1+\lambda_1)(1+\lambda_2)}.
\]
We may solve this with respect to $r$. The other two formulas can be shown in the same way.
\end{proof}
\end{Prop}

\section{Proof of the main theorem}
\setcounter{equation}{0}
In this section, we will show the theorem stated in Introduction.
\setcounter{section}{1}
\begin{Thm}
Assume that a nonsingular curve 
\[
	C: x^4 + y^4 + z^4 + rx^2y^2 + sy^2z^2 + tz^2x^2 = 0
\]
is superspecial, then $r,s$ and $t$ belong to $\mathbb{F}_{p^2}$. Moreover $C$ is maximal or minimal over $\mathbb{F}_{p^2}$.
\end{Thm}\vspace{-1mm}

\setcounter{section}{4}
\setcounter{Def}{0}
A key to the proof is the following proposition by Auer and Top \cite[Proposition 2.2]{Top}.\vspace{-1mm}
\begin{Prop}
Let $E: y^2 = x(x-1)(x-\lambda)$ be a supersingular elliptic curve, then $\lambda \in (\mathbb{F}_{p^2}\hspace{-0.3mm})^4$. Moreover, the followings are true:\vspace{-1mm}
\begin{itemize}
\item If $p \equiv 3 \pmod{4}$, then an elliptic curve $E$ is maximal over $\mathbb{F}_{p^2}$.\vspace{-1mm}
\item If $p \equiv 1 \pmod{4}$, then an elliptic curve $E$ is minimal over $\mathbb{F}_{p^2}$.
\end{itemize}
In particular, the curve $C$ is maximal or minimal over $\mathbb{F}_{p^2}$.
\end{Prop}

\newpage
Recall the discussions in Section 3. Let $E_1,E_2$ and $E_3$ be the following three elliptic curves:
\begin{align*}
	E_1 &: X^2 + y^4 + z^4 + rXy^2 + sy^2z^2 + tz^2X = 0,\\
	E_2 &: x^4 + Y^2 + z^4 + rx^2Y + sYz^2 + tz^2x^2 = 0,\\
	E_3 &: x^4 + y^4 + Z^2 + rx^2y^2 + sy^2Z + tZx^2 = 0.
\end{align*}
Then, there exist surjectives $C \rightarrow E_i$ defined over $\mathbb{F}_{p^2}$.

\begin{proof}[Proof of the first half of Theorem 1.1]
The elliptic curves $E_i$ are supersingular by the assumption, since the quotient of supersingular curve is supersingular. By Lemma 3.1, each elliptic curve $E_i$ is isomorphic to $y^2 = x(x-1)(x-\lambda_i)$. By using Proposition 4.1, each $\lambda_i$ is a fourth power in $(\mathbb{F}_{p^2}\hspace{-0.3mm})^\times$ and thus $\hspace{-0.3mm}\sqrt{\lambda_i} \in \mathbb{F}_{p^2}$. Hence, it follows from Proposition 3.2 that $r,s,t \in \mathbb{F}_{p^2}$.
\end{proof}

Therefore, the question of whether a Ciani curve $C$ is maximal or minimal over $\mathbb{F}_{p^2}$ makes sense. To prove the second assertion of Theorem 1.1, we need the following four lemmas:
\begin{Lem}
We choose $\varDelta$ such that $\varDelta^2 = r^2+s^2+t^2 - rst- 4$, then the followings are true:
\begin{enumerate}
\item If $C$ is superspecial, then $\alpha\beta, \beta\gamma$ and $\gamma\alpha$ belong to $\mathbb{F}_{p^2}$ where $\alpha,\beta,\gamma$ are chosen in (3.1).\vspace{-1mm}
\item If $C$ is superspecial, then $\varDelta$ belongs to $\mathbb{F}_{p^2}$.
\end{enumerate}
\begin{proof}
The elliptic curve $E_i$ is supersingular by assumption, and so $\lambda_i \in (\mathbb{F}_{p^2}\hspace{-0.3mm})^4$ by Proposition 4.1.\par
(1) This claim holds from (3.2) and the first assertion of Theorem 1.1.\par
(2) We obtain $-\lambda_1 \in \mathbb{F}_{p^2}$ clearly, and so whether $\gamma\alpha+(rt-2s)$ is a square in $\mathbb{F}_{p^2}$ is in accord with whether $\gamma\alpha-(rt-2s)$ is a square in $\mathbb{F}_{p^2}$. This means that\vspace{-1.5mm}
\[
	\bigl\{\hspace{-0.3mm}\gamma\alpha+(rt-2s)\bigr\}\bigl\{\hspace{-0.3mm}\gamma\alpha-(rt-2s)\bigr\} = -4(r^2+s^2+t^2-rst-4)\vspace{-1mm}
\]
is a square in $\mathbb{F}_{p^2}$.
\end{proof}
\end{Lem}
\begin{Lem}
If a Ciani curve $C$ is superspecial, then $\frac{r^2-4}{t^2-4},\frac{s^2-4}{r^2-4}$ and $\frac{t^2-4}{s^2-4}$ are fourth powers in $\mathbb{F}_{p^2}$.
\begin{proof}
One can check that
\[
	r^2-4 = \frac{D}{\lambda_1\lambda_2(1-\lambda_3)^2}, \quad s^2-4 = \frac{D}{\lambda_2\lambda_3(1-\lambda_1)^2}, \quad t^2-4 = \frac{D}{\lambda_3\lambda_1(1-\lambda_2)^2}
\]
with $D = (\lambda_1\lambda_2 + \lambda_2\lambda_3 + \lambda_3\lambda_1 - 1)^2 - 4\lambda_1\lambda_2\lambda_3(\lambda_1+\lambda_2+\lambda_3 - 2)$. Hence, we have
\[
	\frac{r^2-4}{t^2-4} = \frac{\lambda_3(1-\lambda_2)^2}{\lambda_2(1-\lambda_3)^2}, \quad \frac{s^2-4}{r^2-4} = \frac{\lambda_1(1-\lambda_3)^2}{\lambda_3(1-\lambda_1)^2}, \quad \frac{t^2-4}{s^2-4} = \frac{\lambda_2(1-\lambda_1)^2}{\lambda_1(1-\lambda_2)^2}.
\]
Here, we have $\lambda_i, 1-\lambda_i \in (\mathbb{F}_{p^2}\hspace{-0.3mm})^4$ by Proposition 4.1, thus $\frac{r^2-4}{t^2-4},\frac{s^2-4}{r^2-4},\frac{t^2-4}{s^2-4}$ are fourth powers in $\mathbb{F}_{p^2}$. Indeed, the elliptic curve $y^2 = x(x-1)(x-\lambda_i)$ is isomorphic to $y^2 = x(x-1)(x-(1-\lambda_i))$ by the proof of \cite[Proposition\,III.1.7]{Silverman}. Therefore $y^2 = x(x-1)(x-(1-\lambda_i))$ is supersingular.
\end{proof}
\end{Lem}
Next, we construct other elliptic curves ${E'}_{\!i}$ which are 2-isogenous to $E_i$. We define
\begin{align*}
	\mu_1 := (rt-2s) + 2\varDelta, \quad \nu_1 := (rt-2s) - 2\varDelta,\\
	\mu_2 := (sr-2t) + 2\varDelta, \quad \nu_2 := (sr-2t) - 2\varDelta,\\
	\hspace{0.2mm}\mu_3 := (ts-2r) + 2\varDelta, \quad \nu_3 := (ts-2r) - 2\varDelta,
\end{align*}
then we obtain $\mu_i,\nu_i \in \mathbb{F}_{p^2\hspace{-0.3mm}}$ by Lemma 4.2\,(2) if a Ciani curve $C$ is nonsingular and superspecial.
\begin{Lem}
Suppose that a Ciani curve $C$ is superspecial. Let ${E'}_{\!i}$ be the elliptic curve defined by\vspace{-0.5mm}
\[
	{E'}_{\!i}: \mu_i\hspace{0.5mm}y^2 = x(x-1)(x-\nu_i/\mu_i).
\]
Then, there exists an isogeny $E_i \rightarrow {E'}_{\!i}$ defined over $\mathbb{F}_{p^2\hspace{-0.3mm}}$ of degree $2$.
\begin{proof}
In this proof, we will show only about ${E'}_{\!1}$. Firstly, the elliptic curve\vspace{-1mm}
\[
	E_1: X^2 + y^4 + z^4 + rXy^2 + sy^2z^2 + tz^2X = 0\vspace{-1mm}
\]
is isomorphic to
\begin{equation}
	v^2 = (r^2-4)u^4 + 2(rt-2s)u^2 + (t^2-4)
\end{equation}
by replacing $u = y$ and $v = 2X+ru^2 + t$. The elliptic curve of the form (4.1) is 2-isogenous to
\begin{equation}
	v^2 = (r^2-4)u^3 + 2(rt-2s)u^2 + (t^2-4)u
\end{equation}
with the morphism $(u,v) \mapsto (u^2,uv)$ defined over $\mathbb{F}_{p^2}$. The right hand side of (4.2) can be factorized as $(r^2-4)u(u-{\mu}'_{1})(u-{\nu}'_{1})$, where
\[
	{\mu}'_{1} = \frac{(rt-2s) + 2\varDelta}{r^2-4}, \quad {\nu}'_{1} = \frac{(rt-2s) - 2\varDelta}{r^2-4}.
\]
Here, the elliptic curve $v^2 = (r^2-4)u(u-{\mu}'_{1})(u-{\nu}'_{1})$ is isomorphic to
\[
	{E'}_{\!1}: \mu_1y^2 = x(x-1)(x-{\nu}'_{1}/{\mu}'_{1}) = x(x-1)(x-\nu_1/\mu_1).
\]
This is the desired conclusion.
\end{proof}
\end{Lem}
There exists an isogeny ${\rm Jac}(C) \rightarrow {E'}_{\!1} \hspace{-0.3mm}\times\hspace{-0.3mm} {E'}_{\!2} \hspace{-0.3mm}\times\hspace{-0.3mm} {E'}_{\!3}$ defined over $\mathbb{F}_{p^2}$, whose degree is a power of $2$. In particular, a curve $C$ is maximal (resp.\ minimal) if and only if ${E'}_{\!i}$ are maximal (resp.\ minimal).

\begin{Lem}
If $C$ is superspecial, then all $\mu_i$ are squares in $\mathbb{F}_{p^2}$ or none of $\mu_i$ is a square in $\mathbb{F}_{p^2}$.
\begin{proof}
It suffices to show that products $\mu_1\mu_2,\mu_2\mu_3,\mu_3\mu_1$ are squares in $\mathbb{F}_{p^2}$. Recall from (3.5) that
\begin{align*}
	\mu_1\nu_1 &= (rt-2s)^2 - 4(r^2+s^2+t^2-rst-4) = (r^2-4)(t^2-4),\\
	\mu_2\nu_2 &= (sr-2t)^2 - 4(r^2+s^2+t^2-rst-4) = (s^2-4)(r^2-4),\\
	\mu_3\nu_3 &= (ts-2r)^2 - 4(r^2+s^2+t^2-rst-4) = (t^2-4)(s^2-4).
\end{align*}
Since ${E'}_{\!i}$ is supersingular by assumption, thus $\nu_i/\mu_i$ is a fourth power in $(\mathbb{F}_{p^2}\hspace{-0.3mm})^\times$ for all $i \in \{1,2,3\}$ by Proposition 4.1. Hence, the product\vspace{-1mm}
\begin{align*}
	\frac{\nu_1}{\mu_1} \cdot \frac{\nu_2}{\mu_2} = \frac{\mu_1\nu_1\mu_2\nu_2}{(\mu_1\mu_2)^2} = \frac{(r^2-4)^2(s^2-4)(t^2-4)}{(\mu_1\mu_2)^2}
\end{align*}
is a fourth power. On the other hand,
\[
	(r^2-4)^2(s^2-4)(t^2-4) = (r^2-4)^4 \cdot \frac{s^2-4}{r^2-4} \cdot \frac{t^2-4}{r^2-4}
\]
is a fourth power by Lemma 4.3, so thus $(\mu_1\mu_2)^2 \in (\mathbb{F}_{p^2}\hspace{-0.3mm})^4$. This means that $\mu_1\mu_2$ is a square in $\mathbb{F}_{p^2}$. In the same way, the products $\mu_2\mu_3$ and $\mu_3\mu_1$ are also squares in $\mathbb{F}_{p^2}$.
\end{proof}
\end{Lem}

Now, we show that a nonsingular superspecial Ciani curve $C$ is maximal or minimal over $\mathbb{F}_{p^2}$.
\begin{proof}[Proof of second half of Theorem 1.1]
Suppose that a Ciani curve $C$ is superspecial, then the elliptic curves ${E'}_{\!1},{E'}_{\!2}$ and ${E'}_{\!3}$ are all supersingular.
\begin{itemize}
\item If $p \equiv 3 \pmod{4}$, then the elliptic curve ${E'}_{\!i}$ is maximal if and only if $\mu_i$ is a square in $\mathbb{F}_{p^2}$, and moreover ${E'}_{\!i}$ is minimal if and only if $\mu_i$ is not a square in $\mathbb{F}_{p^2}$ by using Proposition 4.1. It follows from Proposition 4.5 that all ${E'}_{\!i}$ are maximal or minimal, so the proof is done.\vspace{-1mm}
\item If $p \equiv 1 \pmod{4}$, then the elliptic curve ${E'}_{\!i}$ is minimal if and only if $\mu_i$ is a square in $\mathbb{F}_{p^2}$, and moreover ${E'}_{\!i}$ is maximal if and only if $\mu_i$ is not a square in $\mathbb{F}_{p^2}$ by using Proposition 4.1. It follows from Proposition 4.5 that all ${E'}_{\!i}$ are minimal or maximal, so the proof is done.
\end{itemize}\vspace{-7.7mm}
\end{proof}
\begin{Cor}
Suppose that a Ciani curve $C$ is superspecial, then the followings are true:
\begin{itemize}
\item If $p \equiv 3 \pmod{4}$, then $C$ is maximal if and only if $\mu_i$ for an (equivalently, every) $i \in \{1,2,3\}$ is a square in $\mathbb{F}_{p^2}$.\vspace{-1mm}
\item If $p \equiv 1 \pmod{4}$, then $C$ is maximal if and only if $\mu_i$ for an (equivalently, every) $i \in \{1,2,3\}$ is not a square in $\mathbb{F}_{p^2}$.
\end{itemize}
\end{Cor}
\begin{Exp}
Consider the Ciani curve $C$ defined by the equation
\[
	C: x^4 + y^4 + z^4 + \frac{3}{\sqrt{-2}}x^2y^2 - \frac{9}{4}y^2z^2 + \frac{3}{\sqrt{-2}}z^2x^2 = 0.
\]
A simple computation says $\lambda_1 = -1$ and $\lambda_2 = \lambda_3 = 2$ by using Lemma 3.1. Hence, the three elliptic curves associated to the Ciani curve $C$ are given as below:\vspace{-1mm}
\begin{align*}
	E_1: y^2 = x(x-1)(x+1),\\
	E_2: y^2 = x(x-1)(x-2),\\
	E_3: y^2 = x(x-1)(x-2).
\end{align*}
It follows from \cite[Example\,V.4.5]{Silverman} that these elliptic curves with $j(E_i) = 1728$ are supersingular if and only if $p \equiv 3 \pmod 4$. Therefore $C$ is superspecial if and only if $p \equiv 3 \pmod 4$. Moreover,  the Ciani curve $C$ is maximal over $\mathbb{F}_{p^2}$ if and only if $p \equiv 3 \pmod 4$ since $\mu_1 = \frac{17\sqrt{-1}}{2} \in (\mathbb{F}_{p^2})^2$.
\end{Exp}

Lastly, let us discuss the maximality of a Ciani curve $C$ over $\mathbb{F}_{p^{2e}}$.
\begin{Cor}
Assume that a curve $C: x^4 + y^4 + z^4 + rx^2y^2 + sy^2z^2 + tz^2x^2 = 0$ is superspecial.\vspace{-1mm}
\begin{enumerate}
\item When $e$ is odd, then the followings are true:\vspace{-1mm}
\begin{itemize}
\item If $p \equiv 3 \pmod{4}$, then $C$ is maximal over $\mathbb{F}_{p^{2e}}\!$ if and only if  $\mu_i$ for an (equivalently, every) $i \in \{1,2,3\}$ is a square in $\mathbb{F}_{p^2}$.\vspace{-1mm}
\item If $p \equiv 1 \pmod{4}$, then $C$ is maximal over $\mathbb{F}_{p^{2e}}\!$ if and only if $\mu_i$ for an (equivalently, every) $i \in \{1,2,3\}$ is not a square in $\mathbb{F}_{p^2}$.\vspace{-1mm}
\end{itemize}
\item When $e$ is even, then $C$ is not maximal (i.e. minimal) over $\mathbb{F}_{p^{2e}}$.
\end{enumerate}
\begin{proof}
Let $\varGamma$ be a maximal (resp. minimal) curve over $\mathbb{F}_{p^2}$, then $\varGamma$ over $\mathbb{F}_{p^{2e}}\!$ is also maximal (resp. minimal) if $e$ is odd, and is minimal if $e$ is even. This follows immediately from the Weil conjecture (cf. \cite[Appendix\hspace{1mm}C, Exercise 5.7]{Hartshorne}) and the fact that $\varGamma$ over $\mathbb{F}_{p^2}$ is maximal (resp. minimal) if and only if all the eigenvalues of the Frobenius on the first \'{e}tale cohomology group are $-p$ (resp. $p$).
\end{proof}
\end{Cor}

\end{document}